\documentclass[review]{elsarticle}

\usepackage{lineno,hyperref}
\modulolinenumbers[5]

\journal{Nonlinear Analysis: Hybrid Systems}

\usepackage{amssymb,amsmath,amsthm,color,psfrag}
\usepackage{algorithm,algorithmic,xspace}
\usepackage{subfig}
\usepackage[toc,page]{appendix}

\newcommand{\real}{{\mathbb{R}}}
\newcommand{\norm}[2]{\|#1\|_{#2}}

\newcommand{\subj}{\text{subj. to}}
\newcommand{\half}{\frac{1}{2}}
\newcommand{\ponewt}{\texttt{PRONTO}}

\newcommand{\PP}{\mathcal{P}}
\newcommand{\TT}{\mathcal{T}}

\newtheorem*{remark}{Remark}
\newtheorem{theorem}{Theorem}

\newtheorem{definition}{Definition}

\newcommand\oprocendsymbol{\hbox{$\square$}}
\newcommand\oprocend{\relax\ifmmode\else\unskip\hfill\fi\oprocendsymbol}
    
\usepackage{color}
\newcommand{\Sara}[1]{\textcolor{black}{#1}}%
\newcommand{\SaraBis}[1]{\textcolor{black}{#1}}%

\begin{document}

\begin{frontmatter}

\title{An optimal control approach to the design of periodic orbits for mechanical systems with impacts\tnoteref{erc}}

\author[mymainaddress]{Sara Spedicato\corref{mycorrespondingauthor}}
\ead{sara.spedicato@unisalento.it}

\author[mymainaddress]{Giuseppe Notarstefano}
\ead{giuseppe.notarstefano@unisalento.it}

\cortext[mycorrespondingauthor]{Corresponding author}

\tnotetext[ecr]{This result is part of a project that has received funding from the European Research Council (ERC) under the European Union's Horizon 2020 research and innovation programme (grant agreement No 638992 - OPT4SMART).}

\address[mymainaddress]{Department of Engineering, Universit\`a del Salento, Via per Monteroni, 73100 Lecce, Italy}

\begin{abstract}
In this paper we study the problem of designing periodic orbits for a special class of hybrid systems, namely mechanical systems with underactuated continuous dynamics and impulse events. We approach the
problem by means of optimal control. Specifically, we design an optimal control based strategy that combines trajectory optimization, dynamics embedding, optimal control relaxation and root finding techniques. The proposed strategy allows us to design, in a numerically stable manner, trajectories that optimize a desired cost and satisfy boundary state constraints consistent with a periodic orbit. To show the effectiveness of the proposed strategy, we perform numerical computations on a compass biped model with torso.
\end{abstract}


\end{frontmatter}

\linenumbers

\section{Introduction}

\nolinenumbers

Hybrid systems, involving both continuous and discrete dynamics, arise
naturally in a number of engineering applications.  In particular, many
robotics tasks, such as legged locomotion, (multi-finger) manipulation and
load transportation with aerial vehicles, can be modeled as mechanical systems
with impacts, a particular class of hybrid systems.  These classes of systems
experience continuous dynamics until an interaction with the surrounding
environment (i.e., an impact) occurs, thus causing a state discontinuity due to
impulsive forces and (possibly) a switch between different dynamics.

In this paper we concentrate on a particular aspect that has important
implications both for modeling and control of robots, namely trajectory
design.

\paragraph{Motivations}

The trajectory design task has gained a lot of attention both as a
preliminary task for control and as a tool to explore and understand the
capabilities of the system \Sara{ 
\cite{egerstedt2001optimal},\cite{shaikh2003optimal} ( see \cite{betts1998survey} for an earlier reference).}
Optimization techniques allow us to design
reference trajectories for robot controllers by minimizing a given cost
function. Typical cost functions are (i) the distance from a desired
state-input curve (which does not satisfy the dynamics) and (ii) the energy
injected into the system. 
For example, for humanoid robot design, the distance from a desired (but unfeasible) human-like walking pattern (\cite{vasudevan2013persistent}, \cite{IL-OK-JHP:14}) is often considered.
Additional challenges arise when the trajectory
generation problem is addressed for (underactuated) mechanical systems with
impacts.  The impact events complicate the trajectory optimization problem
since discontinuous changes in the state occur. For some systems, the
underactuation and the instability of the continuous dynamics render the
problem even more challenging.
\Sara{The optimal control theory offers powerfull tools to deal with trajectory generation problems for hybrid dynamical systems.}

\paragraph{Literature review}
The literature on optimal control of hybrid systems is quite vast, thus we report only two sets of contributions relevant for our work.

Fist, a general overview on recent advances for optimal control of hybrid systems is presented.
A recent survey, focusing in particular on switched systems, is
\cite{zhu2015optimal}.
In \cite{baotic2006constrained} the authors propose optimal control algorithms
for discrete-time linear hybrid systems which ``combine a dynamic programming
exploration strategy with multiparametric linear programming and basic
polyhedral manipulation''.  
In \cite{shaikh2007hybrid} a set of necessary conditions is formulated and
optimization algorithms are presented for optimal control of hybrid systems with
continuous nonlinear dynamics and with autonomous and controlled switchings.
The algorithm described in \cite{passenberg2013optimal} (for hybrid systems with
partitioned state space and autonomous switching) is based on a version of the
minimum principle for hybrid systems providing optimality conditions for
intersections and corners of switching manifolds, and thus avoiding the
combinatorial complexity of other algorithms (e.g., \cite{shaikh2007hybrid}).
Furthermore, the work in \cite{wardi2014switched} deals with optimal mode-scheduling
via a gradient descent algorithm for the particular class of autonomous
switched-mode hybrid dynamical systems.
In \cite{riedinger2014lq} and \cite{corona2014stabilization} the design of switching laws for switched systems with linear dynamics, based on the optimization of a quadratic criterion, is addressed.

Second, we focus on trajectory optimization for the particular class of mechanical systems with impacts. 
Besides results regarding systems with impulsive controls,
e.g. \cite{huang2012lq}, \cite{long2011optimal}, we focus on the trajectory design for
systems controlled by inputs of the continuous dynamics.
\Sara{The majority of works adopt parametric optimization methods, as, e.g., \cite{DD-CC-YA:05}, \cite{batz2010ball}, \cite{pchelkin2015algorithms}. 
This means that trajectories are approximated as, for example,  classical, trigonometric or B\'ezier polynomials and the optimization is performed with polynomial coefficients as decision variables. In other works, e.g. \cite{liu2013biped}, the optimization problem is addressed via dynamics equation discretization and the optimal periodic trajectory is computed by means of an approximated cost function. Available software toolboxes are used to solve trajectory generation problems addressed in \cite{DD-CC-YA:05}, \cite{batz2010ball}, \cite{pchelkin2015algorithms} and \cite{liu2013biped}.}
Recently, the optimization over jump times and/or mode sequence, as e.g., in \cite{nakanishi2013spatio} and \cite{posa2014direct}, has been considered.    
In \cite{nakanishi2013spatio} the instants of jump are optimized, but the sequence
of continuous dynamics modes is fixed. In \cite{posa2014direct} the mode
sequence is also optimized by sequential quadratic programming. 
\Sara{
Focusing in particular on the trajectory generation problem for legged robots (the major robotic application regarding mechanical systems with impacts),
challenges arising when dealing with underactuated robots are addressed in  \cite{sugianto2013motion} and \cite{erez2012trajectory} by considering an additional ``virtual'' input acting on the unactuated degrees of freedom. Then, an optimal approximated trajectory of the underactuated dynamics is computed by dynamics inversion.}
Furthermore, online trajectory optimization is addressed in \cite{tassa2012synthesis} through a method based on iterative LQG and in \cite{PLH-AS-LF-UM-SG:13} trajectories are designed by using
an auxiliary system of differential equations and then stabilizing the generated
curves through a feedback on the target system.

\paragraph{Contributions}

As main contribution of this paper, we develop an optimal control based strategy
to design periodic orbits for a class of hybrid dynamical systems with impacts. First, we formulate an optimal control problem with ad-hoc boundary
conditions. These ones are provided by studying how the initial conditions and
the jump conditions are related in a periodic orbit with one jump per
period.  Second, \Sara{instead of using available softwares}, we develop
\Sara{an ad-hoc} strategy based on the combination of trajectory functional
optimization with three main tools: dynamics embedding, constraint relaxation
and zero finding techniques. \Sara{These tools enable us to (i) deal with the
  undeartuated nature of dynamics, (ii) consider highly non trivial constraints
  and (iii) avoid the tedious search for a suitable ``initial (guess)
  trajectory'' to initialize the optimization algorithm.  Furthermore, optimal
  control problems involved in our strategy are solved by combining the penalty
  function approach \cite{bryson1975applied} with the Projection Operator Newton
  method \cite{JH:02}.  In contrast with many strategies reported in the
  literature, we do not resort to approximations such as considering discrete
  sets of motion primitives (thus ending up with the only optimization of their
  parameters), and/or discrete time. In fact we consider system states and
  inputs as optimization variables and a second-order approximation of the
  optimization problem is directly constructed in continuous time.} In detail,
the proposed strategy is based on the following steps. \Sara{By adding a
  fictitious input, we embed the system into a completely controllable (fully
  actuated) one.}  On this system we set up, \Sara{by constraint relaxation}, an
unconstrained optimal control problem to find a trajectory of the system (a
curve satisfying the dynamics) that minimizes a weighted $L_2$ distance from a
desired curve. The desired curve, together with the weights of the cost
functional become important parameters in the designer's hand to explore the system dynamics, that is different periodic orbits that the system
can execute.  \Sara{We solve the optimal control problem to generate a
  trajectory of the fully actuated system with almost null fictitious
  input. Differently from \cite{sugianto2013motion} and
  \cite{erez2012trajectory}, we use this trajectory to initialize the last
  step of our strategy, where the optimal control problem is set up on the
  underactuated dynamics.} A Newton update rule on the final penalty of the
  weighted cost functional is adopted in order to hit the desired final state.
It relies on the solution of an optimal state transfer problem presented in
\cite{JH:03}, where the advantages of the method are highlighted with respect to
the classic ones.

We provide a set of numerical computations showing the effectiveness of the
proposed strategy. In particular, we generate a periodic \Sara{gait} for a three-link
biped robot (compass model with torso). We perform two computations by choosing
two different sets of weights. In the first scenario we try to generate a
trajectory that is as close as possible to the guessed state curve. Vice-versa,
in the second one, we compute a trajectory that minimizes the input effort
(i.e., some sort of minimum-energy trajectory).  
It is worth noting that, although the three links model is relatively simple,
its underactuation represents a significant challenge. Also, it is well known
that for many applications, even such a reduced model of the biped dynamics is
instrumental for analyzing and controlling the actual system, see, e.g.,
\cite{PLH-AS-LF-UM-SG:13}.

\paragraph{Paper organization}
    The paper is organized as follows. In Section~\ref{sec:desired_curve} we introduce 
    the model of the particular class of hybrid system we study in this paper and we present the problem formulation for the generation of periodic orbits. In
    Section~\ref{sec:strategy} we describe the proposed strategy and in
    Section~\ref{sec:example} 
    we provide numerical computations (on a biped walking model) showing its
    effectiveness.

 \section{Problem formulation for the generation of periodic orbits}
 \label{sec:desired_curve}
 In this section we provide the optimal control formulation for solving the problem of periodic orbit generation.

\subsection{Hybrid system model: underactuated systems with impacts}
\label{sec:hybrid_model}
 In this subsection we introduce the particular class of hybrid systems
  studied in the paper.  The hybrid model consists of a continuous dynamics
  phase and a discrete (impulsive) jump event, occurred when the system state
  reaches a jump set.
  \Sara{
  The mechanical system dynamics is modeled by the ordinary
  differential equation
  \begin{equation}
  M(q) \ddot{q} + C(q,\dot{q}) + G(q) = Y_u(q) \; u,
  \label{eq:mech_under}
  \end{equation}
  where $q$, $M(q)$, $C(q,\dot{q})$, $G(q)$, $u$, are respectively the generalized coordinate vector, the mass matrix, the Coriolis vector, the gravity vector, and the input vector.
  The dynamics \eqref{eq:mech_under} can be written in the state space form as
  the continuous and control affine dynamics
    \begin{equation}
    \dot{x}(t)=f(x(t))+g(x(t))u(t),
    \label{eq:state_space}
    \end{equation}
    with state $x \in \real^n$ defined as $x :=[q^T \; \dot{q}^T]^T$ and input $u \in \real^{m}$.  In particular, we
    assume that
    (i) the system is underactuated, i.e., we let $m < n$, and (ii)
    functions $f: \real^n \rightarrow \real^n$ and $g: \real^n \rightarrow \real^{n \times m}$ are (at least) $C^3$.
    }
    When the system state, evolving as modeled in \eqref{eq:state_space},
    reaches a jump set $S \subset \real^{n}$, a discrete impulsive event occurs,
    causing a discontinuity in the system state evolution.
    The discrete event is modeled by the jump map $\Delta: \real^{n} \rightarrow \real^{n}$, which we assume to be invertible, such that
    \begin{equation}
    x^{+}(t) = \Delta (x^{-}(t)),
    \label{eq:discrete}
    \end{equation} 
    where $ x^-(t):=\lim_{\tau \rightarrow t^-} x(\tau) $ and
    $ x^+(t):=\lim_{\tau \rightarrow t^+} x(\tau) $ denote the left and right
    limits of system trajectories satisfying \eqref{eq:state_space}.
    Finally, the overall hybrid model is:
    \begin{equation}
    \Sigma\Sara{:}=\left\{
    \begin{array}{cc}
    \dot{x}(t)=f(x(t))+g(x(t))u(t), & \; x^-(t) \notin S \\
    x^+(t)=\Delta(x^-(t)),      & \; x^-(t) \in S.
    \end{array}
    \right.
    \label{ibrido}
    \end{equation}

\subsection{Optimal control problem set-up} 
We aim at designing periodic trajectories, \Sara{minimizing a weighted $L_2$ distance from a desired curve}, whose period $T$ is fixed. We
approach the task by formulating an optimal control problem over the
continuous dynamics \eqref{eq:state_space}. We impose a boundary
constraint $x(T)=x_f$, where $x_f$ is a final state that allows the system to
jump exactly to the initial condition $x_0$, i.e., $x_0 = \Delta(x_f)$. In
  particular, we choose the initial state $x_0$ and we compute the final state
  by inverting the jump map, i.e., $x_f = \Delta^{-1}(x_0)$.
 
 \Sara{Let a desired curve $(x_d(\cdot),u_d(\cdot))$, satisfying $x_d(0) = x_0$, $x_d(T) = x_f$ and thus $x_d(T) = \Delta^{-1}(x_d(0))$, be given. Furthermore, let us assume that the desired curve \SaraBis{does not reach} the jump set $S$ in the open time interval $(0,T)$.}
 We are ready now to present the optimal control problem we aim to solve for the
 generation of periodic orbits,
 \begin{equation}
 \begin{split}
 \min_{x(\cdot),u(\cdot)} &\;  \half\int_0^T \!\!(\norm{x(\tau) - x_d(\tau)}{Q}^2 +
 \norm{u(\tau) - u_d(\tau)}{R}^2) \; d\tau\\
 \!\!\!\Sara{\subj} &\;  \dot{x}(t) = f(x(t)) + g(x(t)) \; u(t)\\
 & \;  x(0) = x_0\\
 & \;  x(T) = x_f,
 \end{split}
 \label{eq:constr_opt_contr}
 \end{equation}
 \Sara{
where 
$x(\cdot)$ is an absolutely continuous state
trajectory, $u(\cdot)$ is a bounded (measurable) input,  
$Q$ and $R$ are positive definite weighting matrices} and for some vector $z\in\real^k$ and matrix $W\in\real^{k\times k}$ we \Sara{denote
 $\| z \|_W^2 = z^T W z$}.
        
 The time-horizon $T$ is a design parameter together with the desired curve, the
 initial (or final) conditions and the weighting matrices $Q$ and $R$. We stress the fact that the desired
 curve is not a trajectory of the system, i.e., even if it satisfies the
boundary conditions, it does not satisfy the dynamics. Thus, the desired curves
can be seen as ``design tools'' to parameterize actual trajectories of the
system.  \Sara{Furthermore, by choosing the desired trajectory
  $(x_d(\cdot), u_d(\cdot))$ to satisfy the jump condition only at time $T$, we
  expect to compute an optimal (state) trajectory satisfying the same property.
  Otherwise, this constraint can be (implicitly) enforced by properly choosing
  the weight matrices $Q$ and $R$.}
    
\begin{remark}[Final state constraint and controllability]
 It is worth noticing that for problem~\eqref{eq:constr_opt_contr} to be
 feasible the state $x_f$ must belong to the reachable space of the
 system. \oprocend
\end{remark}

\section{Optimal control based strategy for periodic orbit design}
\label{sec:strategy}

\Sara{As highlighted in the literature review, several software toolboxes are
available to solve nonlinear constrained optimal control problems.} Nevertheless, optimal control pro\-blems are optimization problems in the
infinite dimensional space of state-input curves subject to highly non trivial
constraints as the dynamics and the fixed final state constraint. This means
that the solution of an algorithm is highly influenced by the desired curve and
by the ``initial curve'' (possibly a trajectory) that initializes the
algorithm. For this reason, instead of attacking the problem directly by simply
applying one of the solvers available in the literature, we develop a strategy
based on suitable embedding and relaxation ideas that gives a systematic method
to compute periodic trajectories.  

We divide our strategy into three steps that will be presented in detail in the
next subsections. Here we provide an informal idea. First, we add a fictitious
input to the underactuated dynamics, thus getting a fully actuated system. Second, we
consider an \emph{unconstrained} relaxation of the optimal control problem
\eqref{eq:constr_opt_contr} by substituting the final state constraint with a
penalty in the cost function. We solve the relaxed problem by using also the
fictitious input, but with a high penalty. Third and final, we get rid of the
fictitious input and enforce the final state constraint exactly by means of a
suitable iterative method (i.e., a zero finding Newton iteration on the final
state penalty).

\subsection{Dynamics embedding}
The first part of the strategy is called dynamics embedding and is inspired by
\cite{JH-AS-RF:04}. It consists of embedding the underactuated system into a fully
actuated one by adding a fictitious input $u_{emb} \in \real^p$, with $p = n-m$.
\Sara{Thus, embedding the mechanical system \eqref{eq:mech_under} into a fully actuated one, we get
    \begin{equation}
    M(q) \ddot{q} + C(q,\dot{q}) + G(q) = Y(q) \; u^e,
    \label{eq:mech_ful}
    \end{equation}
    where $u^e :=[u^T \; u_{emb}^T]^T$ and $Y(q):= [ Y_u(q) \; Y_{emb}(q)]$ is assumed to be invertible $\forall q$. Accordingly, the state space equations of
    the (embedding) fully-actuated dynamics are
    \begin{equation}
    \dot{x}(t) = f(x(t)) + g(x(t)) u(t) + g_{emb}(x(t)) u_{emb}(t),
    \label{eq:fully_actuated}
    \end{equation}
    where $g_{emb}$ is a \Sara{$\mathcal{C}^3$} vector field taking into account the action of the embedding input.} The role of this additional input is to allow for the tracking of any
(sufficiently regular) desired curve. This provides a useful tool to initialize
the optimal control problem solved in the next part of the strategy.

\subsection{Optimal control problem relaxation and continuation with respect to parameters}
Next, we design an optimal control relaxation of problem
\eqref{eq:constr_opt_contr}. The relaxation involves two aspects. First,
according to the dynamics embedding introduced in the previous step we consider
the fully actuated version of the system. We add a penalty in the cost
function for the embedding input. Second, we relax the final state constraint by
adding a penalty in the cost functional. The relaxed problem is
\begin{equation}
\begin{split}
\min_{x(\cdot),u(\cdot)} &\;  \half\int_0^T (\norm{x(\tau) - x_d(\tau)}{Q}^2 +
\norm{u(\tau) - u_d(\tau)}{R}^2\\
&\; \quad   + \rho_{emb}^2\norm{u_{emb}(\tau)}{}^2) \;  d\tau + \half\rho_{f}^2\norm{x(T) - x_T}{}^2\\
\Sara{\subj} &\;  \dot{x}(t) = f(x(t)) + g (x(t)) u(t) + g_{emb}(x(t)) u_{emb}(t)\\
& \;  x(0) = x_0,
\end{split}
\label{eq:relaxed_opt_contr}
\end{equation}
where $\rho_{emb}$ and $\rho_f$ in $\real_{>0}$ are respectively the weights on
the embedding input and the final state error. The target state $x_T$ is a
strategy parameter that is set to $x_d(T)$ at beginning of the strategy, but
will be modified during the strategy evolution according to the root finding
procedure described in the next step.
The idea is to weight the embedding input and the final state error much
more than the other inputs and the states. In order to avoid a conflict between
the two objectives we propose a continuation procedure on the two penalties.

Problem~\eqref{eq:relaxed_opt_contr} is an unconstrained optimal control
problem that could be solved by means of several available tools in the
literature. We use the \emph{PRojection Operator Newton method for Trajectory Optimization
  (PRONTO)} introduced in \cite{JH:02} and described in \ref{app:pronto}. We adopt the PRONTO tool because it shows two main appealing features for our
trajectory design task. First, it allows one to handle unstable dynamics in a
numerically stable manner and with a ``reasonable'' computational
effort. Second, it guarantees recursive feasibility during the algorithm
evolution. That is, at each step a system trajectory is available. 

We are now ready to present the second step of the strategy. First, given a
desired state curve 
$x_d(\cdot)$ \Sara{such that $x_d =[q_d^T \; \dot{q}_d^T]^T$}
, we can easily compute input trajectories of
the fully actuated system such that the desired state curve is a system
trajectory.  \Sara{By inverting the fully actuated dynamics model
\eqref{eq:mech_ful}, we compute $u_d^e(t)$, $\forall t \in [0,T]$ as
\begin{equation}
u_d^e = Y(q_d)^{-1} [M(q_d) \ddot{q}_d + C(q_d,\dot{q_d}) + G(q_d)].
\nonumber
\end{equation}
Note that, $u_d^e(\cdot)$ \SaraBis{depends on} the choice of $Y_{emb}(\cdot)$ but once $Y_{emb}(\cdot)$ is defined, the input $u_d^e(\cdot)$ is unique.}

Thus, the second step of the strategy can be informally described as follows.
\SaraBis{We denote by} $\xi^e=(x(\cdot),[u^T(\cdot) \; u_{emb}^T(\cdot)]^T)$ a
generic state-input curve for the fully actuated system.  We set the desired
curve to $\xi_d^e=(x_d(\cdot), [u_{d}^T(\cdot) \;
0^T(\cdot)]^T)$ and let \Sara{$\xi^e_0 = (x_d(\cdot), u_d^e(\cdot))$} be the initial trajectory (of the fully
actuated system) for the Projection Operator Newton method (the optimal control
solver). 
Then, \Sara{following an integral penalty function approach \cite{bryson1975applied},} we iteratively solve problem \eqref{eq:relaxed_opt_contr} increasing the weight
$\rho_{emb}$ at each step, \Sara{by means of a suitable heuristic}. When the norm of the embedding input is sufficiently small 
\Sara{(i.e., when $\norm{u_{emb}(\cdot)}{} <\epsilon_{emb}$, where $\epsilon_{emb}$ is a given tolerance)} we stop
the procedure, \Sara{thus obtaining an approximated trajectory of the underactuated system}.

A pseudo-code description of the strategy is given in the next table (Algorithm 1).
We denote $\PP^{e}$ the projection operator \eqref{eq:proj_oper} acting on the fully actuated system,
so that $\eta^e = \PP^{e}(\xi^e)$ is a trajectory of the fully-actuated
system. For given $\rho_{emb}$ and $\rho_f$, we denote $\ponewt^{e}$ a routine such that
$\xi^e_{opt^e} = \ponewt^{e}(\xi^e_0, \xi^e_d; \rho_{emb}, \rho_f)$, i.e., it  
takes as inputs
an initial trajectory $\xi^e_0$ and a desired (state-control) curve $\xi^e_d$,
and computes an embedding optimal trajectory, which is denoted as $\xi^e_{opt^e}=(x_{opt^e}(\cdot),
[u_{opt^e}^T(\cdot) \; u_{opt^e,emb}^T(\cdot)]^T)$, solving the nonlinear
optimal control relaxation in \eqref{eq:relaxed_opt_contr}. 

\begin{algorithm}[H]
    \caption{Periodic orbit design strategy: step 2}
    \label{alg:gait_design_strategy_2}
    \begin{algorithmic}
        \REQUIRE Initial condition $x_0$ \Sara{and embedding desired curve $\xi_d^{e}$}
        \STATE compute initial trajectory $\xi_{0}^{e} = \PP^{e}(\xi_{0}^{e})$
        \FOR{$i=0, 1, \ldots$}
        \STATE \% solve relaxed optimal control for given $\rho_{emb}$\\
        $\xi^{e}_{opt^e,i}$ = $\ponewt^{e}(\xi_i^{e},
        \xi_{d}^{e}; \rho_{emb}, \rho_f$)\\
        \STATE {\bf if} $\norm{u_{emb}(\cdot)}{} <\epsilon_{emb}$  exit for  {\bf end if}\\
        \STATE increase $\rho_{emb}$
        \STATE set $\xi_{i+1}^{e}=\xi^{e}_{opt^e,i}$
        \ENDFOR
        \ENSURE $ \xi^{e}_{opt^e} = \xi^{e}_{opt^e,i}$
    \end{algorithmic}
\end{algorithm}
\begin{remark}
\SaraBis{
(Convergence)
Provided that a solution to \eqref{eq:relaxed_opt_contr} with $u_{emb}(\cdot)=0$ exists, the existence of the solution is also guaranteed, by continuity, when $u_{emb}$ is in a neighbourhood of the origin (with suitable radius $r_{emb} \geq \epsilon_{emb}$). Thus, a solution to \eqref{eq:relaxed_opt_contr} with $||u_{emb}(\cdot)||<\epsilon_{emb}$ exists and the convergence to that solution follows since we are using an integral penalty function approach (see, e.g. \cite{bryson1975applied}). Nevertheless, dealing with a non-convex problem, only a local convergence is ensured, provided that the algorithm is initialized with a trajectory lying inside the basin of attraction.}
\end{remark}

\subsection{Enforcing the final state constraint}
We are ready to present the third and final step of the strategy. This part
relies on the idea, proposed in \cite{JH:03}, of enforcing the final state
constraint of problem \eqref{eq:constr_opt_contr} by combining two actions: (i)
increase the penalty $\rho_f$ on the final state error, and (ii) vary the target
state $x_T$ until the candidate optimal trajectory ``practically'' meets the
final state constraint.

The following optimal control problem is the relaxation
of problem \eqref{eq:constr_opt_contr} subject to the original system,
i.e., without dynamics embedding
\begin{equation}
\begin{split}
\min_{x(\cdot),u(\cdot)} &\;  \half\int_0^T \!\!\!(\norm{x(\tau) - x_d(\tau)}{Q}^2 +
\norm{u(\tau) - u_d(\tau)}{R}^2) \; d\tau  \\
& + \half\rho_{f}^2\norm{x(T) - x_T}{}^2\\
\!\!\Sara{\subj} &\;  \dot{x}(t) = f(x(t)) + g (x(t)) u(t)\\
& \;  x(0) = x_0.
\end{split}
\label{eq:relaxed_opt_contr2}
\end{equation}

\Sara{Before presenting the theorem that establishes a connection between the
  original problem \eqref{eq:constr_opt_contr} and its relaxed version
  \eqref{eq:relaxed_opt_contr2}, let us state the second order sufficient
  condition for optimality satisfied by a local minimum of the original problem
  \eqref{eq:constr_opt_contr}.
\begin{definition} \textit{(Second order sufficient condition \cite{JH:03})}
The second order sufficient condition for a local minimum is
\begin{equation*}
\int_0^T
\left[
\begin{array}{c}
z(\tau)\\
v(\tau)\\
\end{array}
\right]^T
\left[
\begin{array}{cc}
H_{xx}(\tau) & H_{xu}(\tau)\\
H_{xu}^T(\tau) & H_{uu}(\tau)\\
\end{array}
\right] 
\left[
\begin{array}{c}
z(\tau)\\
v(\tau)\\
\end{array}
\right]
d\tau \; > 0, \quad \forall z \in \real^n, \forall v \in \real^m,
\end{equation*}
where $H$ is the control Hamiltonian and $H_{uu} = \frac{\partial^2 H}{\partial u^2}$, $H_{xx}= \frac{\partial^2 H}{\partial x^2}$,  $H_{xu} = \frac{\partial^2 H}{\partial x \partial u}$.
\end{definition}}
The following result is at the basis of the third step of
the proposed strategy.

\begin{theorem}\textit{(Equivalence of constrained and relaxed minimizers \cite{JH:03})}
    Suppose that $\xi^* = (x^*(\cdot), u^*(\cdot))$ is a local minimum of
    \eqref{eq:constr_opt_contr} satisfying the second order sufficiency condition
    for optimality. Then for each $\rho_f>0$, there is a target state $x_T =
    x_T(\rho_f)$ such that $\xi^*$ is a local minimizer of problem
    \eqref{eq:relaxed_opt_contr2} satisfying the
    corresponding second order sufficiency condition. \oprocend
\end{theorem}

Thus, the third step of the strategy can be informally described as
follows. Given an embedding optimal trajectory $\xi^e_{opt^e}$ from the previous
step, we obtain the initial trajectory $\xi_0=\mathcal{P}((x_{opt^e}(\cdot), u_{opt^e}(\cdot))$ of the original system,
where
$\PP$ is the
projection operator acting on the original system.  Furthermore, we reduce the embedding desired curve $\xi_d^e$, in order to get the desired curve  for the original system as 
$\xi_d=(x_d(\cdot),u_{d}(\cdot))$.
Then, \Sara{adopting a penalty function approach,}
we iteratively solve problem
\eqref{eq:relaxed_opt_contr2} increasing the weight $\rho_f$ at each step. 
\Sara{The use of a large $\rho_{f}$ can easily result in a poorly conditioned problem where the terminal cost overwhelms the remaining one. Thus, we start with a reasonably small $\rho_{f}$ and we gradually increase it by means of a suitable heuristic.
Note that our penalty function approach generates an optimal trajectory which approximately satisfies the final state constraint. For this reason,
when the final state $x_{opt,i}(T)$ of the temporary optimal
trajectory is ``sufficiently close'' to the desired final state (i.e., $\norm{x_{opt,i}(T)-x_d(T)}{} \leq \delta^f_{tol}$, where $\delta^f_{tol}$ is tolerance guaranteeing the root finding
convergence), we apply a
Newton method for root finding on $x_T$ to meet exactly the final state
constraint.} \Sara{The update rule we use on $x_T$ is the one proposed in \cite{JH:03}.    
According to the Implicit Function Theorem, (i) 
there exists a mapping $\beta : \real^n \rightarrow \real^n$ such that $x_{opt,i}(T) = \beta(x_{T,i})$, where $\xi_{opt,i}(\cdot) = (x_{opt,i}(\cdot), u_{opt,i}(\cdot))$ is the solution to \eqref{eq:relaxed_opt_contr2} with target state $x_{T,i}$ and (ii) the first Fr\'echet differential of $\beta(\cdot)$ at $x_{T,i}$, i.e., $z_f \mapsto D\beta(x_{T,i}) \cdot z_f$, exists.
Furthermore, provided that the linearization of $\dot{x}(t) = f(x(t)) + g(x(t)) u(t)$ about $\xi_{opt,i}$ is controllable on $[0,T]$, the mapping $\beta(\cdot)$ is invertible. Thus, a Newton method for root finding is applied on the final state constraint equation $\beta(x_T) - x_f = 0$ in order to find the value of $x_T$ such that the solution to \eqref{eq:relaxed_opt_contr2} is equivalent to the solution to \eqref{eq:constr_opt_contr}.
In particular, at the iteration $i+1$,
we update the value of $x_T$ according to the following rule: $x_{T,i+1} = x_{T,i} + D\beta(x_{T,i})^{-1} (x_f - \beta(x_{T,i}))$.
The final state constraint is met when $\norm{x_{opt,i}(T)-x_d(T)}{} <\epsilon^f_{tol}$, where $\epsilon_{tol}^f$ is the desired tolerance on the final
state error.}

A pseudo code description of the third step is given in the
following table (Algorithm 2). We use the following notation.
For a given $\rho_f$ \Sara{and $x_T$}, we denote $\ponewt$ the Projection Operator Newton method
for the original system, so that we have \Sara{$\xi_{opt} = \ponewt(\xi_0, \xi_d;
\rho_f, x_T)$}.
\Sara{
Note that, as for Algorithm 1, convergence of Algorithm 2 directly follows since it is a properly formulated penalty function approach. 
}

\begin{algorithm}[H]
    \caption{Periodic orbit design strategy: step 3}
    \label{alg:gait_design_strategy_3}
    \begin{algorithmic}
        \REQUIRE $\xi^{e}_{opt^e}$ and $\xi^{e}_d$
        \STATE reduce embedding desired curve\\
        \STATE $\xi_d=(x_d(\cdot), u_{d}(\cdot))$\\[1.1ex]
        \STATE \% project $(x_{opt^e}(\cdot), u_{opt^e}(\cdot))$ to get an initial trajectory for $\ponewt$
        \STATE $\xi_0 = \PP((x_{opt^e}(\cdot), u_{opt^e}(\cdot)))$
        \FOR{$i=0, 1, \ldots$}
        \STATE compute: $\xi_{opt,i}$ = $\ponewt(\xi_i,
        \xi_{d}; \rho_{f}, \Sara{x_{T,i}})$\\
        \STATE {\bf if} $\norm{x_{opt,i}(T)-x_d(T)}{} <\epsilon^f_{tol}$  {\bf then} exit  {\bf end if}\\[1.1ex]
        \IF{$\norm{x_{opt,i}(T)-x_d(T)}{} >\delta^f_{tol}$}
        \STATE increase $\rho_{f}$
        \ELSE \STATE \Sara{$x_{T,i+1} = x_{T,i} + D\beta(x_{T,i})^{-1} (x_d(T) - \beta(x_{T,i}))$}
        \ENDIF
        \STATE set: $\xi_{i+1}=\xi_{opt,i}$
        \ENDFOR
        \ENSURE $ \xi_{opt} = \xi_{opt,i}$
    \end{algorithmic}
\end{algorithm}    
        
\section{An illustrative example: three link planar biped robot}
\label{sec:example}

As an example, we apply the proposed technique to a planar biped walking model. The model adopted in this paper was introduced in \cite{JWG-GA-FP:01}.

\subsection{Biped model}

The three rigid links biped with four lumped masses
consists of a torso and two legs of equal length connected at the hip.
The hybrid walking model consists of a continuous swing phase model and a discrete jump event model.

The swing phase model describes the motion of the swing leg to develop a walking step.
Let $\theta_1, \; \theta_2, \; \theta_3$ denote the orientation of the stance leg, the swing leg and the torso, respectively. Furthermore, let $u_1$ be the torque applied between the stance leg and the torso while the torque $u_2$ acts between the swing leg and the torso. 
The swing phase model via Euler-Lagrange equations is given by
\begin{equation}
M(\theta)\ddot{\theta}+C(\theta,\dot{\theta})+G(\theta)=U
\label{eq:constr}
\end{equation}
where $\theta=[\theta_1 \; \theta_2 \; \theta_3]^T$ and $M(\theta)$, $C(\theta,\dot{\theta})$, $G(\theta)$ and $U$ (defined in \ref{app:model_matrices}) are respectively the mass
matrix, the Coriolis vector, the gravity vector and the generalized torques vector.
Defining the state
$
x=[\theta_1 \; \theta_2 \; \theta_3 \; \dot{\theta}_1 \; \dot{\theta}_2 \;
\dot{\theta}_3]^T,
$
and the input $u=[u_1 \; u_2]^T$,
we can write the dynamics \eqref{eq:constr} in state-space form \eqref{eq:state_space}.

The jump event model takes into account (i) the impulse force arising when the swing leg touches the ground and (ii) the switch of the leg being in contact with the ground. See \cite{JWG-GA-FP:01} for more details.
The jump event model is given by \eqref{eq:discrete} where
\begin{equation}
\Delta([\theta^{-T} \; \dot{\theta}^{-T}])=
\left[
\begin{array}{cc}
R &  0_{3 \times 3}\\
0_{3 \times 3} &  A(\theta^-)\\
\end{array}
\right]
\left[
\begin{array}{c}
\theta^-\\
\dot{\theta}^-\\
\end{array}
\right] 
\end{equation}
with $A(\theta^-)$ reported in \ref{app:impact_model} and
\begin{equation*}
R=
\left[
\begin{array}{ccc}
0 &  1 & 0 \\
1 &  0 & 0 \\
0 &  0 & 1 \\
\end{array}
\right].
\end{equation*}
The jump event occurs when the system state reaches the jump set
\begin{equation}
S:= \{x \in \real^6 |
\theta_1=\theta_1^{jmp}\}
\end{equation}
where $\theta_1^{jmp}$ is set according to physical considerations.

\subsection{Numerical computations}
\label{sec:num_comp}
According to \cite{JWG-GA-FP:01}, we consider the
following model parameters: mass of the legs $m=5$ kg; mass of the hip $M_H=15$
kg; mass of the torso $M_T=10$ kg; length of the legs $r=1$ m, and the length of
the torso $l=0.5$ m. We use as jump angle $\theta_1^{jmp} = \pi / 8$.
We choose as time horizon $T=1.53$ and as initial condition $x_0=[- 22.5\;\text{deg} \;\; 22.5\;\text{deg} \;\; 20\;\text{deg} 
\;\; 50\;\text{deg/s} \;\; 0\;\text{deg/s} \;\; 
90\;\text{deg/s}]^T$.

The first step of the proposed strategy requires the definition of a fully actuated dynamics and the design of a desired curve.
In order to obtain a fully actuated dynamics, we
consider the dynamics \eqref{eq:constr} and add a (generalized) torque acting
directly on the torso angle velocity $\dot{\theta}_3$. That is, the swing phase
model becomes
\begin{equation}
M(\theta)\ddot{\theta}+C(\theta,\dot{\theta})+G(\theta)=U^{e},
\label{eq:constr2}
\end{equation}
where $U^e=Y u^e$, 
\begin{equation}
Y=
\left[
\begin{array}{ccc}
-1 &  0 & 0 \\
0 &  -1 & 0 \\
1 &  1 & 1 \\
\end{array}
\right],
\end{equation}
and $u^e=[u_1 \; u_2 \; u_{emb}]^T$.
To design a desired curve that satisfies constraints on the initial and
final states, we choose as desired angles $\theta_{d,k}(\cdot)$ with
$k=1,2,3$ spline curves so that the initial and final point and
their relative slopes, i.e., $\theta_{d,k}(0)$, $\dot{\theta}_{d,k}(0)$,
$\theta_{d,k}(T)$ and $\dot{\theta}_{d,k}(T)$, can be assigned a priori.
Then, the desired velocity and acceleration curves, respectively
$\dot{\theta}_{d}(\cdot)$ and $\ddot{\theta}_{d}(\cdot)$, are obtained by
\emph{symbolic} time differentiation.
According to the
fully-actuated dynamic model \eqref{eq:constr2}, we compute $ \forall t \in
[0,T]$
\[
u_d^e(t)=Y^{-1}(M(\theta_{d}(t))\ddot{\theta}_{d}(t)+C(\theta_{d}(t),\dot{\theta}_{d}(t))+G(\theta_{d}(t))),
\]
where \Sara{$u_d^e=[u_{0}^T  \; u_{0,emb}]^T$, with $u_0 \in \real^2$ and $u_{0,emb} \in \real$. 
Furthermore, as regards the update rule for the penalty parameters, at iteration $i+1$, we set $\rho_{emb,i+1} = 2 \; \rho_{emb,i}$. The same update rule is adopted for $\rho_f$.} Numerical computations characterized by two different design objectives are presented in the following. We also invite the reader to watch
the attached video showing optimal walking gaits obtained by means of the proposed strategy.

\begin{figure*}[htbp!]
    \begin{center}
        \subfloat[][Snapshots of the optimal trajectory (Gait 1)]{
            \includegraphics[width=115mm]{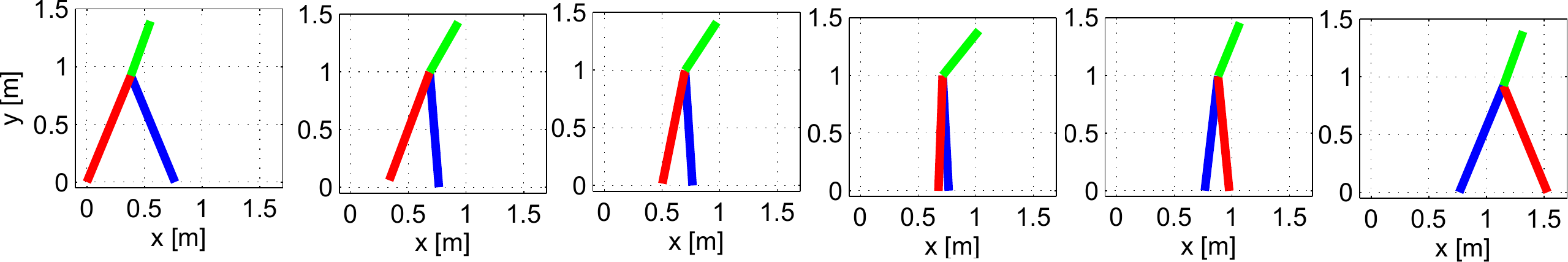}}\\
        \subfloat[][Snapshots of the optimal trajectory (Gait 2)]{
            \includegraphics[width=115mm]{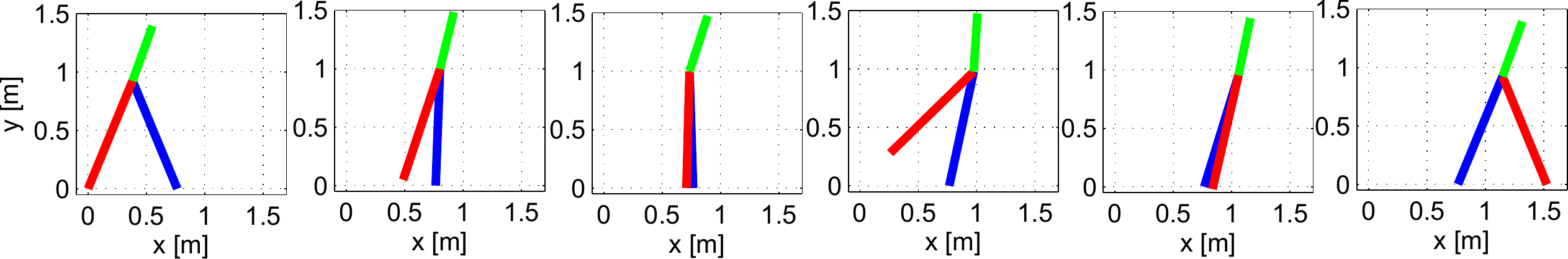}}\quad
        \subfloat[][Stance leg position $\theta_1$ (Gait 1)]{
            \includegraphics[width=50mm]{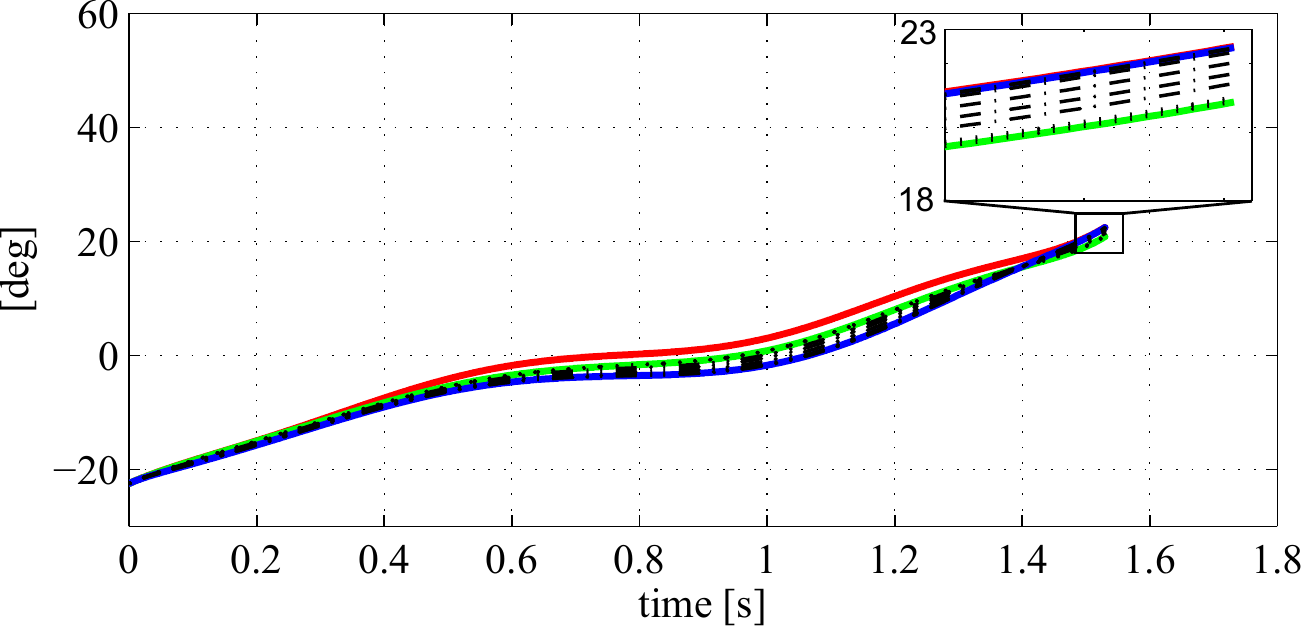}
        } \quad \quad
        \subfloat[][Stance leg position $\theta_1$ (Gait 2)]{
            \includegraphics[width=50mm]{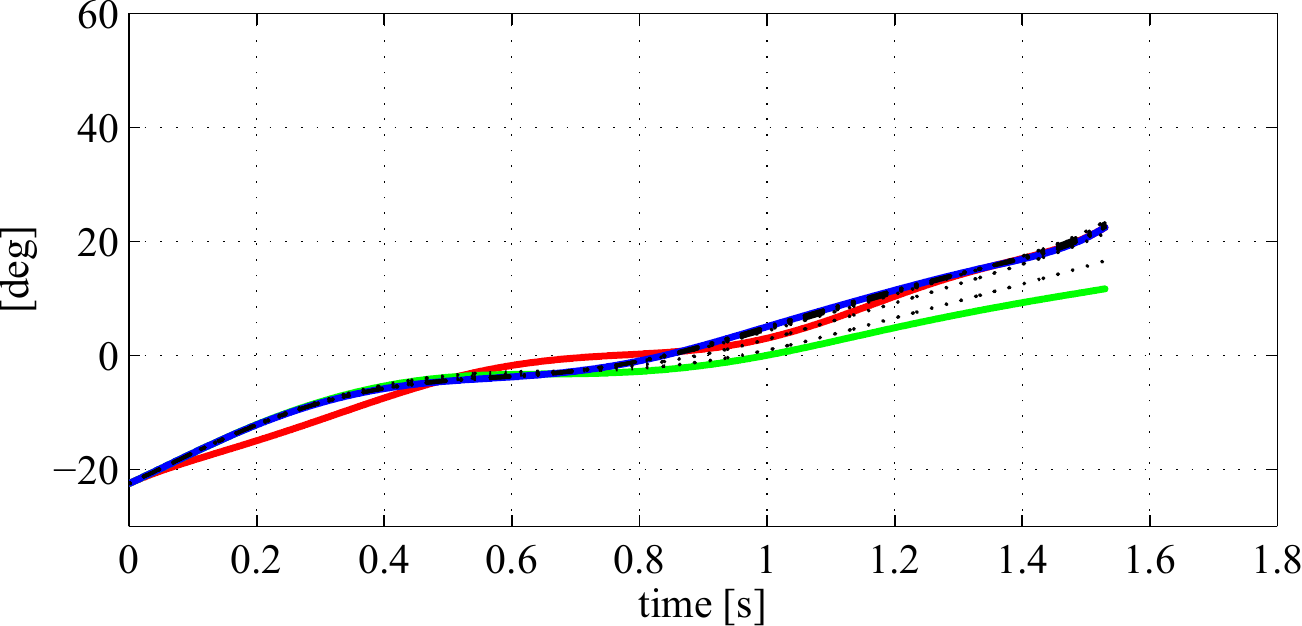}
        }\\
        \subfloat[][Swing leg position $\theta_2$ (Gait 1)]{
            \includegraphics[width=50mm]{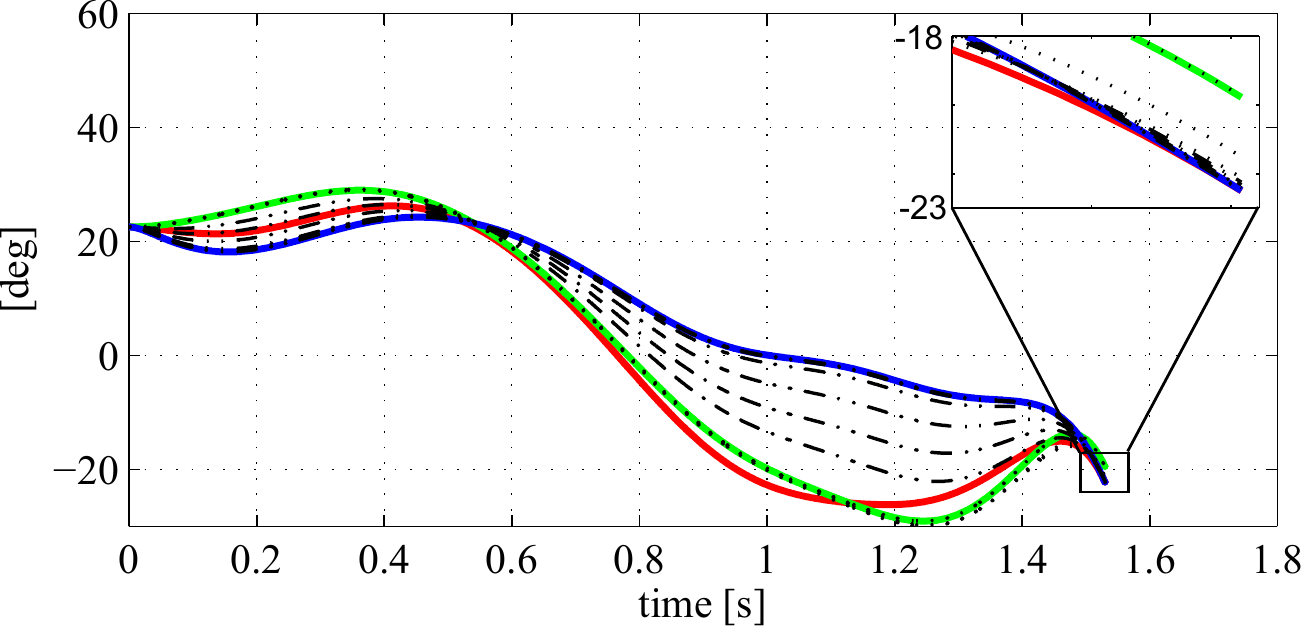}
        } \quad \quad    
        \subfloat[][Swing leg position $\theta_2$ (Gait 2)]{
            \includegraphics[width=50mm]{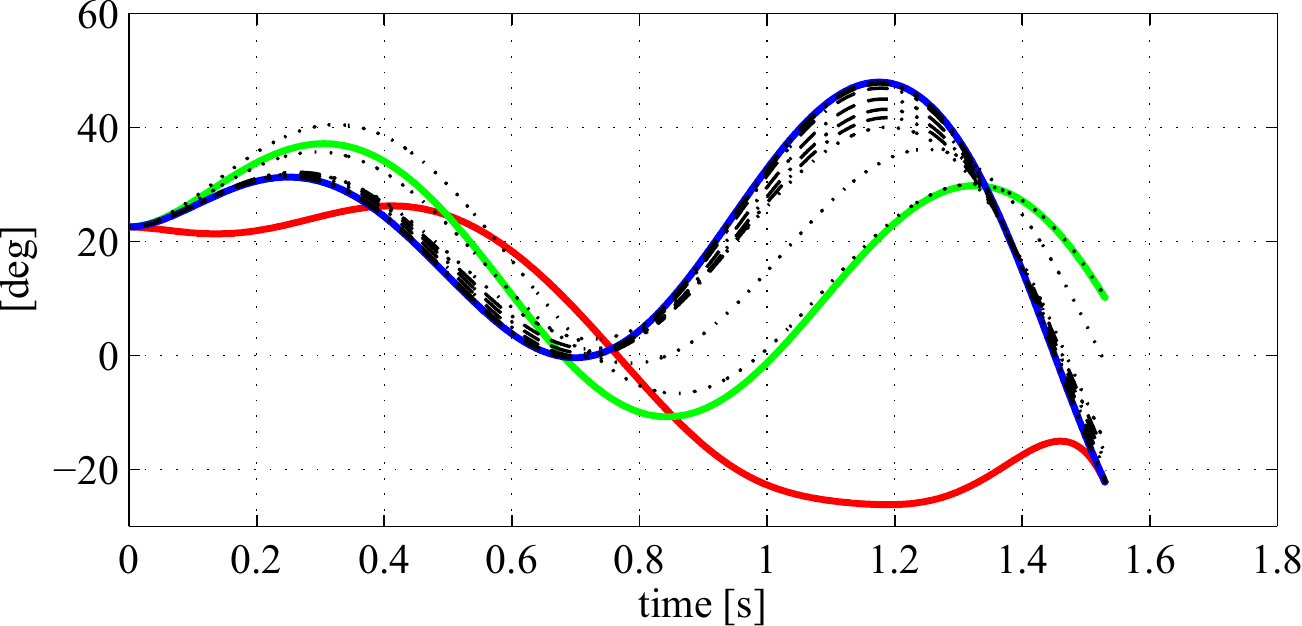}
        }    \\
        \subfloat[][Torso position $\theta_3$ (Gait 1)]{
            \includegraphics[width=50mm]{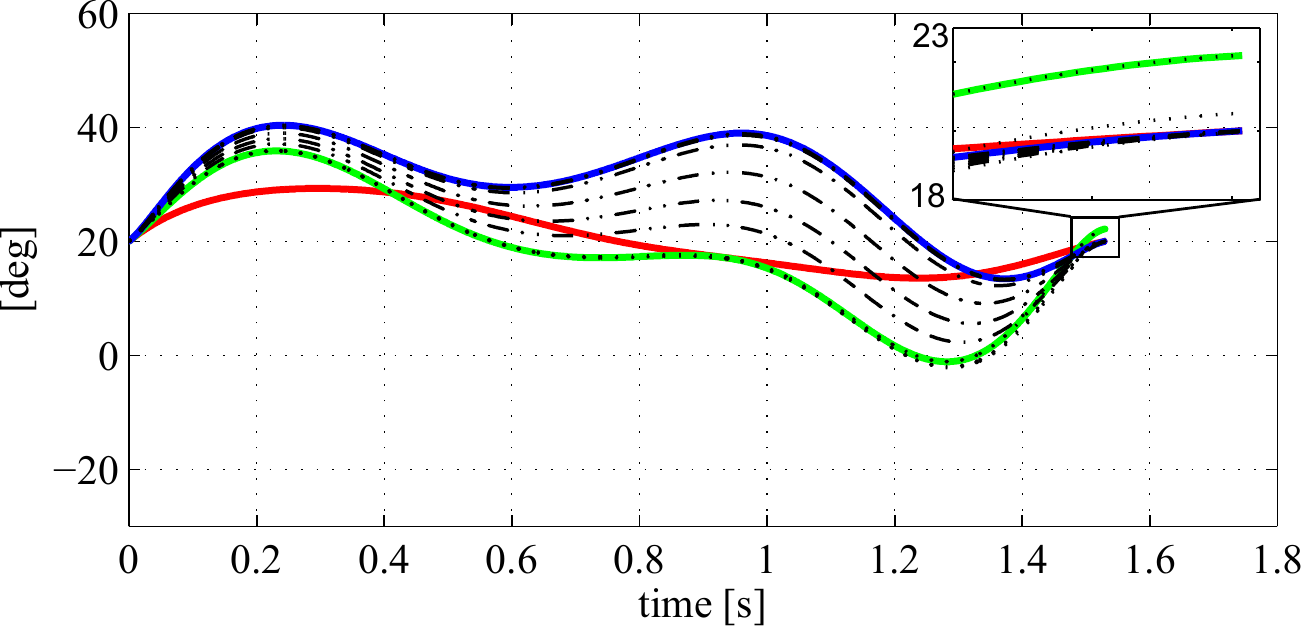}
        }  \quad \quad  
        \subfloat[][Torso position $\theta_3$ (Gait 2)]{
            \includegraphics[width=50mm]{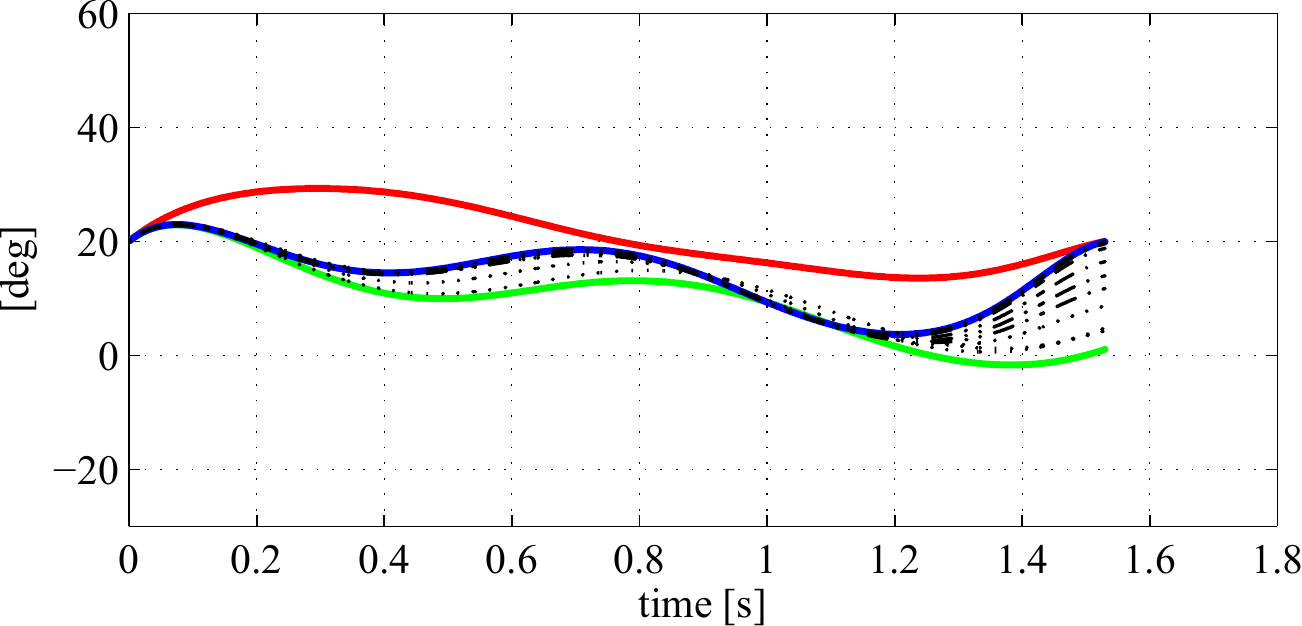}
        }                  
        \caption{Gait 1 (left column) and gait 2 (right column): optimal trajectory (solid blue) compared to the desired infeasible curve (solid red) and the temporary optimal trajectory of the embedded system \Sara{(solid green)}. \Sara{Temporary optimal trajectories of the uderactuated system are depicted in black: trajectories varying $\rho_f$ (dotted line) and updating the target state $x_T$ (dot-dashed line).}}
        \label{fig:states}  
    \end{center}           
\end{figure*}

\begin{figure*}[ht!]
    \begin{center}
        \subfloat[][Input $u_1$ (Gait 1)]{
            \includegraphics[width=50mm]{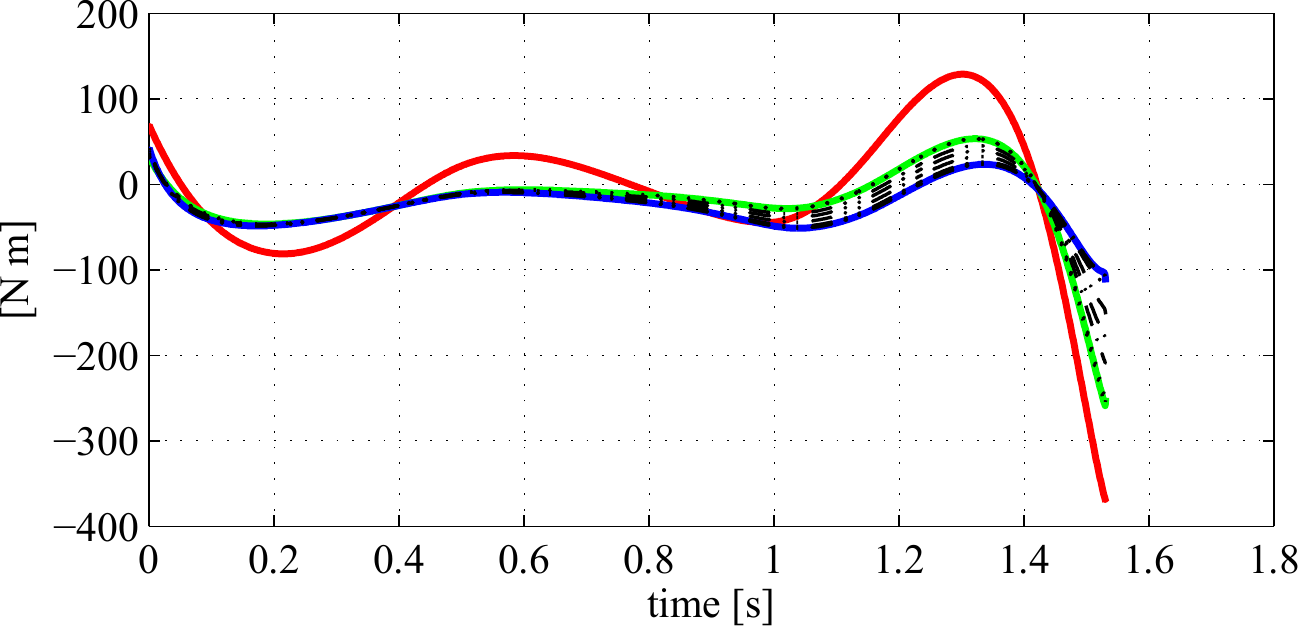}
        } \quad \quad
        \subfloat[][Input $u_1$ (Gait 2)]{
            \includegraphics[width=50mm]{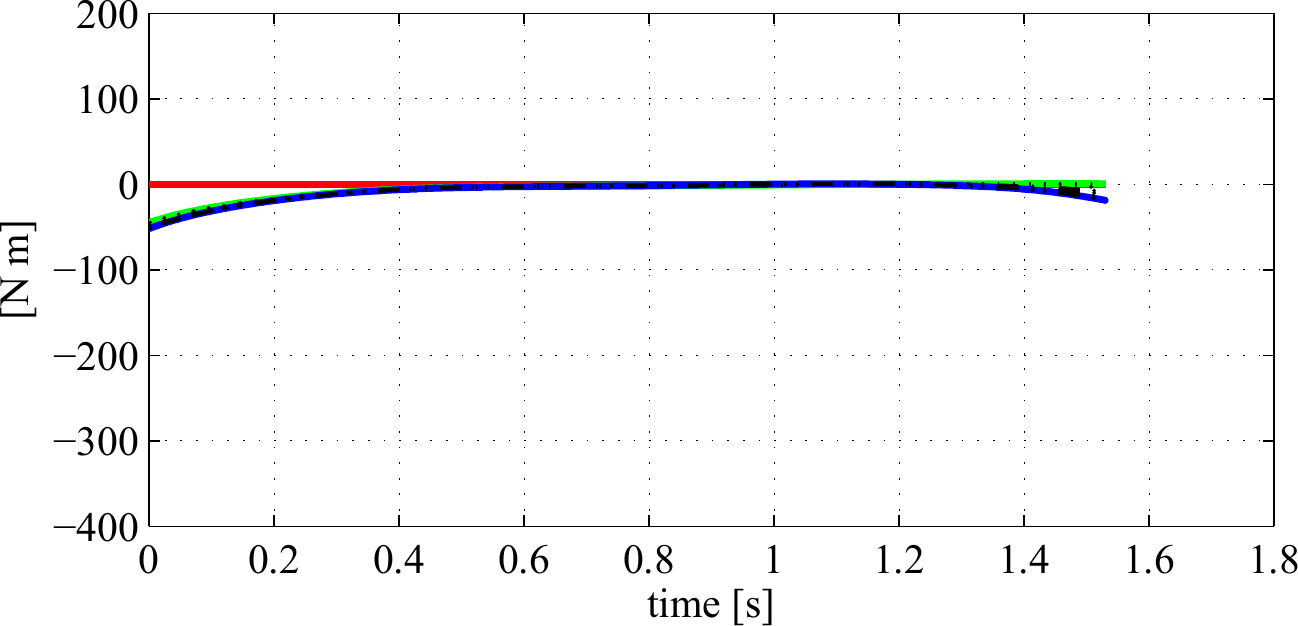}
        }  \\
        \subfloat[][Input $u_2$ (Gait 1)]{
            \includegraphics[width=50mm]{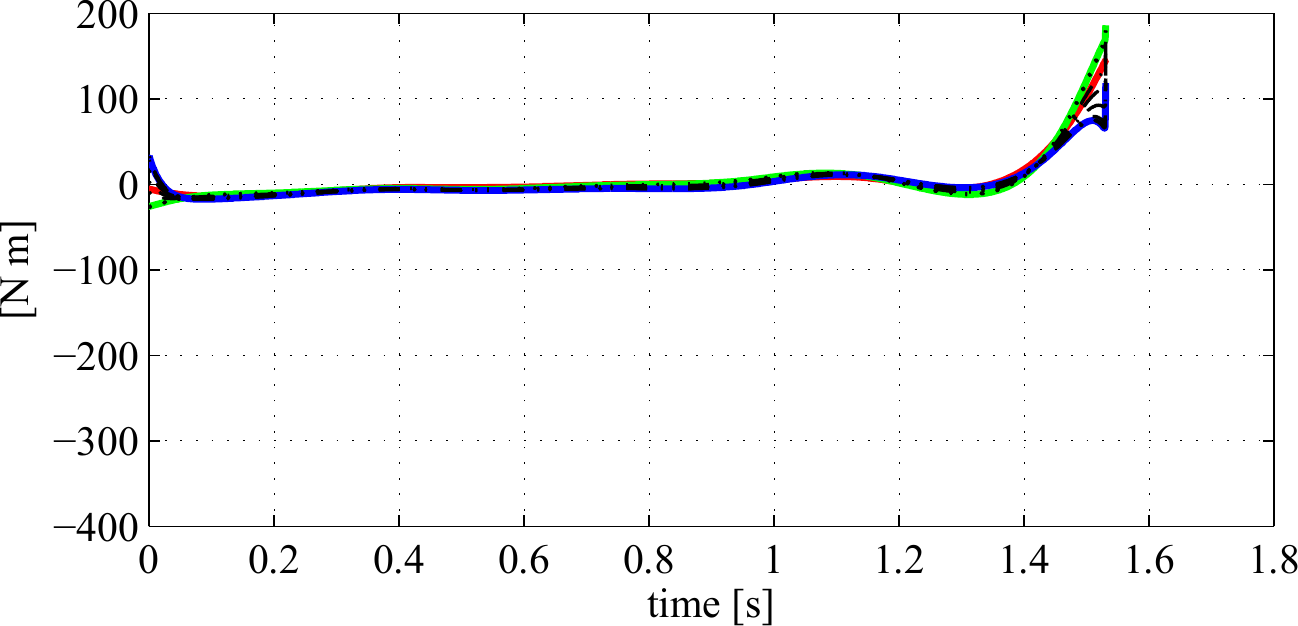}\label{fig:u2_gait1}  
        } \quad \quad
        \subfloat[][Input $u_2$ (Gait 2)]{
            \includegraphics[width=50mm]{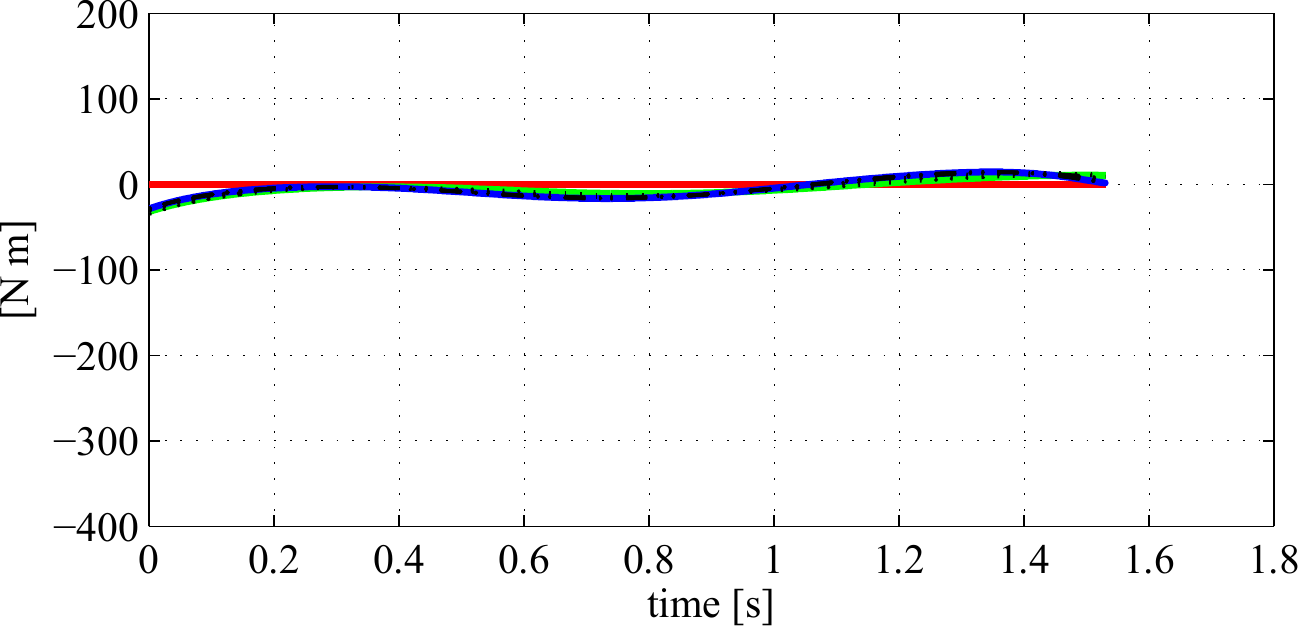}\label{fig:u2_gait2}  
        }                                      
        \caption{Inputs of gait 1 (left column) and gait 2 (right column): optimal trajectory (solid blue) compared to the desired infeasible curve (solid red) and the temporary optimal trajectory of the embedded system \Sara{(solid green)}. \Sara{Temporary optimal trajectories of the underactuated system are depicted in black: trajectories varying $\rho_f$ (dotted line) and updating the target state $x_T$ (dot-dashed line).}}
        \label{fig:inputs}  
    \end{center}           
\end{figure*}

In the first computation, \Sara{we choose $x_d(t) = [\theta_d(t)^T \; \dot{\theta}_d(t)^T ]^T$ and $u_d(t) = u_0(t)$, $\forall t \in [0,T]$}.
We find the optimal trajectory choosing diagonal Q and
R matrices and penalizing the angles $10^4$ times the input, while the
velocities $10^3$ times the input. That is, $Q= \text{diag}
\left[
100 \; 100 \; 100 \; 10 \; 10 \; 10
\right]
$
and
$R=\text{diag}
\left[
0.01 \;  0.01
\right]
$.
The optimal trajectory (blue) is depicted on the left column of
Figure~\ref{fig:states} and Figure~\ref{fig:inputs} 
together with the desired curve (red), \Sara{the embedding optimal trajectory (green) and temporary trajectories of the underactuated system (black).} 
The optimization strategy
shapes the desired curve in order to obtain a trajectory (a curve
satisfying the dynamics), thus showing that in this case the desired curve is
just a rough guess of the biped gait.

In the second computation, we look for a trajectory that minimizes the ``energy'' injected in the system. Indeed, we choose 
\Sara{$x_d(t) = [\theta_d(t)^T \; \dot{\theta}_d(t)^T ]^T$, $u_d(t) = [0 \; 0]^T$, $\forall t \in [0,T]$,}
and penalize
the inputs $10^3$ times the angles and $10^2$ times the velocities. That is, we
choose $Q=\text{diag} \left[
0.01 \;  0.01 \; 0.01 \; 0.1 \; 0.1 \; 0.1
\right]$
and
$R=\text{diag}
\left[
10 \;  10
\right]$.
The optimal trajectory (blue) is depicted 
on the right column of Figure~\ref{fig:states} 
and
Figure~\ref{fig:inputs}
together with the desired curve (red), \Sara{the embedding optimal trajectory (green) and temporary trajectories of the underactuated system (black).}
As expected, the desired curve and the
optimal trajectory coincide at $T$, while the embedding optimal trajectory has a nonzero
final state error.  This shows the effectiveness of the third step of the
strategy that enforces the final state constraint.
This optimal trajectory is
different from the previous one, thus showing another possible gait.
In order to minimize the injected energy, the swing leg first goes behind the stance one, it gains potential energy and then, it passes in front of the stance leg. This gait requires a less amount of torque, with respect to the first gait. Approaching the final time $T$, the torque $u_2$ of the first gait (Figure \ref{fig:u2_gait1}) is greater than $100 \;  N m$ while it approaches zero in the second gait (Figure \ref{fig:u2_gait2}).     

\section{Conclusion}
\label{sec:conclusion}
In this paper we have developed an optimal control based strategy to compute
periodic orbits for underactuated mechanical systems with impacts. The strategy combines trajectory
optimization with dynamics embedding, optimal control relaxation and root
finding techniques. The proposed strategy provides a systematic and numerically
robust methodology to design periodic orbits for a particular class of hybrid systems.

\section*{References}

\bibliography{bibliography}

\newpage

\appendix

\section{The Projection Operator approach for the
    optimization of trajectory functionals}
\label{app:pronto}

\Sara{
The PRojection Operator based Newton method for Trajectory Optimization (PRONTO) \cite{JH:02} is suitable for solving optimal control problems in the form
\begin{equation}
\begin{split}
& \min  \int_0^T l(x(\tau),u(\tau)) \; d\tau + m(x(T))\\
& \subj \quad \dot{x}(t) = f(x(t),u(t)), \quad x(0)=x_0,
\end{split}
\label{eq:PRONTO_problem}
\end{equation}
where the initial condition $x_0$ is fixed and $l(x,u)$, $m(x)$ and $f(x,u)$ are (at least) $\mathcal C^3$ in $x$ and $u$. Sufficient conditions on $f$, $l$ and $m$ (\cite{cesari1983optimization}, \cite{lee1989optimalcontrol}) guarantee the existence of optimal trajectories. In order to deal with state-input constraints, a strategy combining PRONTO with a barrier function approach is proposed in \cite{hauser2006barrier}.  
The key idea of PRONTO is that a properly designed \emph{projection operator} $\mathcal{P}$, mapping  state-control curves into system trajectories (curves satisfying the dynamics), 
is used to convert the dynamically constrained optimization problem \eqref{eq:PRONTO_problem}
into an essentially unconstrained one. 
Let $\xi = (\alpha(\cdot), \mu(\cdot))$ be a bounded curve and let $\eta = (x(\cdot), u(\cdot))$ be a trajectory of the nonlinear feedback system
\begin{equation}
\begin{split}
\dot{x}(t) =& f(x(t), u(t)), \qquad x(0) = x_0,\\
u(t) =& \mu(t) + K(t)(\alpha(t) - x(t)),
\end{split}
\label{eq:proj_oper_def}
\end{equation} 
where the initial condition $x_0$ is given in \eqref{eq:PRONTO_problem} and the feedback gain $K(\cdot)$ is designed, e.g., by solving a suitable linear quadratic optimal control problem, in order to guarantee (local) exponential stability of the trajectory $\eta$. 
The feedback system \eqref{eq:proj_oper_def} defines the nonlinear \emph{projection operator} 
\begin{equation}
\PP : \xi \mapsto \eta,
\label{eq:proj_oper}
\end{equation}
mapping the curve $\xi$ to the trajectory $\eta$.
Using the projection operator to locally parameterize the trajectory manifold, problem \eqref{eq:PRONTO_problem} 
is equivalent to the one of minimizing the unconstrained functional
$g(\xi) = h(\mathcal{P}(\xi))$, where $h(\xi) := \int_0^T l(x(\tau),u(\tau)) \; d\tau + m(x(T))$.
Then, using an (infinite dimensional) Newton descent method, a local minimizer is computed.  The strength of this approach
is that the local minimizer is obtained as the limit of a sequence of
trajectories, i.e., curves satisfying the dynamics.
The feedback system \eqref{eq:proj_oper_def}, defining the projection operator, allows us to generate trajectories in a numerically stable manner. In other words, the choice of the feedback in \eqref{eq:proj_oper_def} is convenient from a numerical point of view. Furthermore, note that (projected) trajectories $(x_i(\cdot), u_i(\cdot))$ satisfy $x_i(0) = x_0$, according to the definition of the projection operator.}

\Sara{
A pseudo-code of the Projection Operator Newton method is shown in the table (Algorithm 3). Let $\xi_0$ be an initial
trajectory. Minimization of the cost functional $g(\xi)$ is accomplished iteratively. 
Given the current trajectory iterate $\xi_i$, the search direction $\zeta_i$ is obtained by solving a
linear quadratic optimal control problem with cost $Dg(\xi_i) \cdot \zeta + \frac{1}{2} D^2 g(\xi_i)(\zeta,\zeta)$, where
$\zeta \mapsto Dg(\xi_i) \cdot \zeta$ and $\zeta \mapsto
D^2 g(\xi_i)(\zeta,\zeta)$ are respectively the first and second Fr\'echet
differentials of the functional $g(\xi) \!\!= \!\!h(\PP(\xi))$ at $\xi_i$. Then,
the curve $\xi_i + \gamma_i \zeta_i$, where $\gamma_i$ is a step size obtained through
a standard backtracking line search, is projected, by means of the projection operator, in order to get a new trajectory $\xi_{i+1}$.
}
\begin{algorithm}[H]
    \caption{Projection Operator Newton method}
    \label{alg:proj_newt}
    \begin{algorithmic}
        \REQUIRE initial trajectory $\xi_0 \in \TT$\\
        
        \FOR{$i = 0, 1, 2 \ldots$}
        \STATE design $K$ defining $\PP$ about $\xi_i$
        \STATE search for descent direction\\
        $\zeta_i = \text{arg} \min_{\zeta\in T_{\xi_i} \TT} Dg(\xi_i) \cdot \zeta + \frac{1}{2} D^2 g(\xi_i)(\zeta,\zeta)$
        \STATE step size $\gamma_i = \arg \min_{\gamma \in (0,1]} g(\xi_i + \gamma \zeta_i)$;\\
        \STATE project $\xi_{i+1}={\PP}(\xi_i + \gamma_i \zeta_i)$.
        \ENDFOR
    \end{algorithmic}
\end{algorithm}
\begin{remark}
    The algorithm has the structure of a standard Newton method for the
    minimization of an unconstrained function. The key points are the design of
    $K$ defining the projection operator and the computation of the derivatives of
    $g$ to ``search for descent direction''. It is worth noting that these two
    steps involve the solution of suitable (well known) linear quadratic optimal
    control problems, \cite{JH:02}.\oprocend
\end{remark}

\section{Biped walking model: swing phase}
\label{app:model_matrices}

    The mass
    matrix, the Coriolis vector, the gravity vector and the generalized torques vector of equation \eqref{eq:constr} are respectively
    \begin{equation}
    M(\theta)= \left[
    \begin{array}{ccc}
    (\frac{5}{4}m+M_H+M_T)r^2 &
    -\frac{1}{2}mr^2c_{12} &
    M_Trlc_{13}\\
    -\frac{1}{2}mr^2c_{12} &
    \frac{1}{4}mr^2 &
    0\\
    M_Trlc_{13} &
    0 &
    M_Tl^2
    \end{array}
    \right]
    \label{M}
    \end{equation}
    \begin{equation}
    C(\theta,\dot{\theta})= \left[
    \begin{array}{c}
    -\frac{1}{2}mr^2s_{12}\dot{\theta_2}^2+ M_Trls_{13}\dot{\theta_3}^2\\
    \frac{1}{2}mr^2s_{12}\dot{\theta_1}^2\\
    -M_Trls_{13}\dot{\theta_1}^2
    \end{array}
    \right]
    \label{eq:C}
    \end{equation}
    \begin{equation}
    G(\theta)= \left[
    \begin{array}{c}
    -\frac{1}{2}g(2M_H+3m+2M_T)r\sin\theta_1 \\
    \frac{1}{2}gmr\sin\theta_2 \\
    -gM_Tl\sin\theta_3
    \end{array}
    \right]
    \label{eq:G}
    \end{equation}
    \begin{equation}
    U= \left[
    \begin{array}{c}
    -u_1 \\
    -u_2 \\
    u_1+u_2
    \end{array}
    \right]
    \label{eq:U}
    \end{equation}
    where
    \begin{equation}
    \begin{array}{c}
    c_{12}:=\cos(\theta_1-\theta_2) \quad c_{13}:=\cos(\theta_1-\theta_3)\\
    s_{12}:=\sin(\theta_1-\theta_2) \quad s_{13}:=\sin(\theta_1-\theta_3),
    \end{array}
    \end{equation}
    and the model parameters are:
    the mass of the legs $m$, the mass of the hip $M_H$, the mass of the torso $M_T$, the length of the legs $r$, and the length of
    the torso $l$.

\section{Biped walking model: relation between velocities\\ just after and before the impact}
\label{app:impact_model}

The relation between velocities just after and before the impact is represented by $\dot{\theta}^+=A(\theta^-)\dot{\theta}^-$. The terms of the matrix
$A(\theta^-)$ are reported below.
$$
A_{11}=\frac{1}{\text{den}}[2M_T\cos(-\theta^{-}_{1}-\theta^{-}_{2}+2\theta^{-}_{3})+
$$
$$
-(2m+4M_H+2M_T)\cos(\theta^{-}_{1}-\theta^{-}_{2}),
$$
$$
A_{12}=\frac{1}{\text{den}}m, \quad A_{13}=0,
$$
$$
A_{21}=\frac{1}{\text{den}}[m-(4m+4M_H+2M_T)\cos(2\theta^{-}_{1}-2\theta^{-}_{2})+
$$
$$
2M_T\cos(2\theta^{-}_{1}-2\theta^{-}_{3})],
$$
$$
A_{22}=\frac{1}{\text{den}}[2m\cos(\theta^{-}_{1}-\theta^{-}_{2})], \quad A_{23}=0,
$$
$$
A_{31}=\frac{1}{l\text{den}}[(2mr+2M_Hr+2M_Tr)\cos(\theta^{-}_{1}-2\theta^{-}_{2}+\theta^{-}_{3})+
$$
$$
-2M_Hr\cos(-\theta^{-}_{1}+\theta^{-}_{3})- (2mr+2M_Tr)
$$
$$
\cos(-\theta^{-}_{1}+\theta^{-}_{3})+mr\cos(-3\theta^{-}_{1}+2\theta^{-}_{2}+\theta^{-}_{3}),
$$
$$
A_{32}=-\frac{1}{l\text{den}}rm\cos(-\theta^{-}_{2}+\theta^{-}_{3}), \quad A_{33}=1,
$$
$$
\text{den}=-3m-4M_H-2M_T+2m\cos(2\theta_{1}^--2\theta_{2}^-)+$$
$$2M_T\cos(-2\theta_{2}^-+2\theta_{3}^-).
$$    
    
\end{document}